\documentclass{amsart}
\usepackage{amssymb}
\usepackage{amsmath}
\usepackage{amsfonts}

\newtheorem{theorem}{Theorem}[section]
\theoremstyle{plain}

\newtheorem{corollary}[theorem]{Corollary}

\newtheorem{lemma}[theorem]{Lemma}

\newtheorem{proposition}[theorem]{Proposition}
\newtheorem{remark}[theorem]{Remark}

\numberwithin{equation}{section}
\input{tcilatex}

\begin{document}
\title[On Schatten-class properties of pseudo-differential operators]{On
Schatten-von Neumann class properties of pseudo-differential operators.
Cordes' lemma.}
\author{Gruia Arsu}
\address{Institute of Mathematics of The Romanian Academy\\
P.O. Box 1-174\\
RO-70700\\
Bucharest \\
Romania}
\email{agmil@home.ro}
\subjclass[2000]{Primary 35S05, 43Axx, 46-XX, 47-XX; Secondary 42B15, 42B35.}
\maketitle

\begin{abstract}
We prove an extended version of Cordes' lemma concerning trace-class
properties of some special pseudo-differential operators. This version of
Cordes' lemma is used to improve the results in \cite{Arsu} concerning the
Schatten-class properties of pseudo-differential operators in the $\left(
X,\tau \right) $-quantization. Here $X$ is an $n$ dimensional vector space
and $\tau $ is an endomorphism of $X$.
\end{abstract}

\section{Introduction}

In \cite{Arsu}, we extend a method due to H.O. Cordes and T. Kato to deal
with Schatten-class properties of pseudo-differential operators. Let $\left(
X,\left\vert \cdot \right\vert _{X}\right) $ be an $n$ dimensional euclidean
space and $\left( X^{\ast },\left\vert \cdot \right\vert _{X^{\ast }}\right) 
$ its dual. We prove, among others, that if a symbol $a\left( x,\xi \right) $
defined on $X\times X^{\ast }$ has $L^{p}$ derivatives $D_{x}^{\alpha
}D_{\xi }^{\beta }a$\ for $\left\vert \alpha \right\vert ,\left\vert \beta
\right\vert \leq $\ $\left[ n/2\right] +1$ and $1\leq p<\infty $, then the
associated pseudo-differential operator $a\left( x,D\right) $ belongs to $%
\mathcal{B}_{p}\left( L^{2}\left( X\right) \right) $ (the Schatten ideal of
compact operators whose singular values lie in $l^{p}$). The result is
actually finer, the conditions being imposed on the derivatives
corresponding to an orthogonal decomposition of $X$.

The extension of Cordes-Kato method can be used to obtain similar results
for any other quantization $a^{\tau }$, where $\tau $ is an endomorphism of
the vector space $X$, if an appropriate $\tau $-version of Cordes' lemma can
be proved.

The purpose of the present paper is to prove an extension of the $\tau $%
-version of Cordes' lemma \ to more general class of symbols and to more
general parameter $\tau $. For example, $\tau $ will belong to an open
neighborhood $\mathcal{U}$ of \ $%
%TCIMACRO{\U{211d} }%
%BeginExpansion
\mathbb{R}
%EndExpansion
\cdot 1_{X}$ in the space of endomorphisms of $X$,\ while the class of
symbols, which contains the special symbols used in the original Cordes'
lemma, will be defined in the next section. The open neighborhood $\mathcal{U%
}$ of \ $%
%TCIMACRO{\U{211d} }%
%BeginExpansion
\mathbb{R}
%EndExpansion
\cdot 1_{X}$ in $\limfunc{End}_{%
%TCIMACRO{\U{211d} }%
%BeginExpansion
\mathbb{R}
%EndExpansion
}\left( X\right) $, the space of endomorphisms of $X$, is defined by $%
\mathcal{U}=\mathcal{U}_{0}\cup \mathcal{U}_{1}$, where $\mathcal{U}_{0}$ is
the set of all invetible endomorphisms of $X$ and $\mathcal{U}_{1}=1_{X}+%
\mathcal{U}_{0}$.

Let us note that the proof in \cite{Cordes} of this lemma was done only in a
particular case. The reason for doing this was to avoid considerable
difficulties of technical nature, due to the complicate structure of the
singular convolution operator $\left( 1-\triangle \right) ^{t}$, $t>0$. It
is of some interest to give a new proof, which is applicable to the general
case. The techniques we use are similar to those used in the appendix of 
\cite{Arsu}.

Finally, this $\tau $-version of Cordes' lemma allow us to improve the
results in \cite{Arsu} concerning the Schatten-class properties of
pseudo-differential operators in the $\left( X,\tau \right) $-quantization.
Here $X$ is an $n$ dimensional vector space and $\tau $ is an endomorphism
of $X$.

\section{Cordes' lemma}

Let $\left( X,\left\vert \cdot \right\vert \right) $ be an euclidean space
of dimension $n$ and $\left( X^{\ast },\left\vert \cdot \right\vert
_{X^{\ast }}\right) $ its dual. If $x\in X$, we set $\left\langle
x\right\rangle =\left( 1+\left\vert x\right\vert ^{2}\right) ^{1/2}$.
Sometimes, in order to avoid confusions, we shall add a subscript specifying
the space, e.g. $\left( \cdot ,\cdot \right) _{X}$, $\left\vert \cdot
\right\vert _{X}$ or $\left\langle \cdot \right\rangle _{X}$.

To state and prove Cordes' lemma we shall work with a very restricted class
of symbols. We shall say that $a:X\rightarrow 
%TCIMACRO{\U{2102} }%
%BeginExpansion
\mathbb{C}
%EndExpansion
$ is a symbol of degree $m$ ($m$ any real number) if $a\in \mathcal{C}%
^{\infty }\left( X\right) $ and for any $k\in 
%TCIMACRO{\U{2115} }%
%BeginExpansion
\mathbb{N}
%EndExpansion
$, there is $C_{k}>0$ such that 
\begin{equation*}
\left\Vert a^{\left( k\right) }\left( x\right) \right\Vert \leq
C_{k}\left\langle x\right\rangle ^{m-k},\quad x\in X.
\end{equation*}%
If an orthonormal basis is given in $X$, this is equivalent with the
requirement that 
\begin{equation*}
\left\vert \partial ^{\alpha }a\left( x\right) \right\vert \leq C_{\alpha
}\left\langle x\right\rangle ^{m-\left\vert \alpha \right\vert },
\end{equation*}%
for all $\alpha \in 
%TCIMACRO{\U{2115} }%
%BeginExpansion
\mathbb{N}
%EndExpansion
^{n}$ and all $x\in X$. We denote by $\mathcal{S}^{m}\left( X\right) $ the
vector space of all symbols of degree $m$ and observe that 
\begin{equation*}
m_{1}\leq m_{2}\Rightarrow \mathcal{S}^{m_{1}}\left( X\right) \subset 
\mathcal{S}^{m_{2}}\left( X\right) ,\quad \mathcal{S}^{m_{1}}\left( X\right)
\cdot \mathcal{S}^{m_{2}}\left( X\right) \subset \mathcal{S}%
^{m_{1}+m_{2}}\left( X\right) .
\end{equation*}%
Observe also that $a\in \mathcal{S}^{m}\left( X\right) \Rightarrow \partial
^{\alpha }a\in \mathcal{S}^{m-\left\vert \alpha \right\vert }\left( X\right) 
$ for each $\alpha \in 
%TCIMACRO{\U{2115} }%
%BeginExpansion
\mathbb{N}
%EndExpansion
^{n}$. The function $\left\langle x\right\rangle ^{m}$ clearly belongs to $%
\mathcal{S}^{m}\left( X\right) $ for any $m\in 
%TCIMACRO{\U{211d} }%
%BeginExpansion
\mathbb{R}
%EndExpansion
$. We denote by $\mathcal{S}^{\infty }\left( X\right) $ the union of all the
spaces $\mathcal{S}^{m}\left( X\right) $ and we note that $\mathcal{S}\left(
X\right) =\bigcap_{m\in 
%TCIMACRO{\U{211d} }%
%BeginExpansion
\mathbb{R}
%EndExpansion
}\mathcal{S}^{m}\left( X\right) $ the space of tempered test functions. It
is clear that $\mathcal{S}^{m}\left( X\right) $ is a Fr\'{e}chet space with
the seni-norms given by 
\begin{equation*}
\left\vert a\right\vert _{m,\alpha }=\sup_{x\in X}\left\langle
x\right\rangle ^{-m+\left\vert \alpha \right\vert }\left\vert \partial
^{\alpha }a\left( x\right) \right\vert ,\quad a\in \mathcal{S}^{m}\left(
X\right) .
\end{equation*}

\begin{lemma}
Let $a\in \mathcal{S}^{0}\left( X\right) $ and set $a_{\varepsilon }\left(
x\right) =a\left( \varepsilon x\right) $. Then $\left\{ a_{\varepsilon
}\right\} _{0\leq \varepsilon \leq 1}$ is bounded in $\mathcal{S}^{0}\left(
X\right) $ and for every $m>0$, $a_{\varepsilon }\rightarrow a_{0}$ in $%
\mathcal{S}^{m}\left( X\right) $ when $\varepsilon \rightarrow 0$.
\end{lemma}

\begin{proof}
The statement follows if we show that for $0<m\leq 1$ 
\begin{equation*}
\left\langle x\right\rangle ^{-m+\left\vert \alpha \right\vert }\left\vert
\partial ^{\alpha }\left( a_{\varepsilon }\left( x\right) -a\left( 0\right)
\right) \right\vert \leq C_{\alpha }\varepsilon ^{m},\quad x\in X,0\leq
\varepsilon \leq 1.
\end{equation*}%
When $\alpha =0$ this follows by Taylor's formula. We have 
\begin{eqnarray*}
\left\vert a\left( x\right) -a\left( 0\right) \right\vert &=&\left\vert
\left\langle x,\int_{0}^{1}a^{\prime }\left( \lambda x\right) \limfunc{d}%
\lambda \right\rangle \right\vert \leq C_{0}\int_{0}^{\left\vert
x\right\vert }\left\langle \lambda \right\rangle ^{-1}\limfunc{d}\lambda \\
&\leq &2^{1/2}C_{0}\ln \left( 1+\left\vert x\right\vert \right) \leq
C_{m}\left\vert x\right\vert ^{m},\quad x\in X,
\end{eqnarray*}%
which implies 
\begin{equation*}
\left\vert a_{\varepsilon }\left( x\right) -a\left( 0\right) \right\vert
=\left\vert a\left( \varepsilon x\right) -a\left( 0\right) \right\vert \leq
C_{m}\varepsilon ^{m}\left\langle x\right\rangle ^{m},\quad x\in X.
\end{equation*}%
When $\alpha \neq 0$ we just have to use that 
\begin{equation*}
\left\langle x\right\rangle ^{-m+\left\vert \alpha \right\vert }\left\langle
\varepsilon x\right\rangle ^{-\left\vert \alpha \right\vert }\varepsilon
^{-m+\left\vert \alpha \right\vert }\leq \left( \left( 1+\left\vert
x\right\vert ^{2}\right) /\left( \varepsilon ^{-2}+\left\vert x\right\vert
^{2}\right) \right) ^{\left( -m+\left\vert \alpha \right\vert \right)
/2}\leq 1.
\end{equation*}
\end{proof}

\begin{corollary}
\label{c3}Let $r\in 
%TCIMACRO{\U{211d} }%
%BeginExpansion
\mathbb{R}
%EndExpansion
$. Then $\bigcup_{\rho <r}\mathcal{S}^{\rho }\left( X\right) \subset 
\overline{\mathcal{S}\left( X\right) }^{\mathcal{S}^{r}\left( X\right) }$,
where $\overline{\mathcal{S}\left( X\right) }^{\mathcal{S}^{r}\left(
X\right) }$ is the closure of $\mathcal{S}\left( X\right) $ in $\mathcal{S}%
^{r}\left( X\right) $.
\end{corollary}

\begin{proof}
Let $\rho <r$ and $b\in \mathcal{S}^{\rho }\left( X\right) $. Choose $\chi
\in \mathcal{C}_{0}^{\infty }\left( X\right) $, $\chi \left( 0\right) =1$.
Then $\left\{ \chi _{\varepsilon }\right\} _{0<\varepsilon \leq 1}$ is
bounded in $\mathcal{S}^{0}\left( X\right) $ and $\chi _{\varepsilon
}\rightarrow \chi _{0}\equiv 1$ in $\mathcal{S}^{r-\rho }\left( X\right) $
when $\varepsilon \rightarrow 0$.\ It follows that $b^{\varepsilon }=\chi
_{\varepsilon }b\in \mathcal{C}_{0}^{\infty }\left( X\right) $ and $%
b^{\varepsilon }-b=\left( \chi _{\varepsilon }-\chi _{0}\right) b\rightarrow
0$ in $\mathcal{S}^{r}\left( X\right) $ when $\varepsilon \rightarrow 0$.
\end{proof}

Let $\varphi \in \mathcal{C}_{0}^{\infty }\left( X^{\ast }\right) $ such
that $0\leq \varphi \leq 1,$ $\varphi \left( p\right) =1$ for $\left\vert
p\right\vert \leq 1$, $\varphi \left( p\right) =0$ for $\left\vert
p\right\vert \geq 2$. Then we have 
\begin{equation*}
\limfunc{d}\varphi \left( p/t\right) /\limfunc{d}t=\psi \left( p/t\right)
/t,\quad \psi \left( p\right) =-\left\langle \nabla \varphi \left( p\right)
,p\right\rangle =-\sum_{j=1}^{n}p_{j}\partial _{j}\varphi \left( p\right)
,\quad p\in X^{\ast },
\end{equation*}%
which yields a continuous partition of unity 
\begin{equation*}
1=\varphi \left( p\right) +\int_{1}^{\infty }\psi \left( \frac{p}{t}\right) 
\frac{\func{d}t}{t},\quad p\in X^{\ast }.
\end{equation*}%
Note that $t\leq \left\vert p\right\vert \leq 2t$ in the support of $\psi
\left( p/t\right) $ ($\limfunc{supp}\psi \subset \left\{ p\in X^{\ast
}:1\leq \left\vert p\right\vert \leq 2\right\} $).

We shal make use of the following simple but important remark. If $a\in 
\mathcal{S}^{m}\left( X^{\ast }\right) $, then the family $\left\{ \psi
a_{t}\right\} _{1\leq t<\infty }$ is bounded in $\mathcal{S}\left( X^{\ast
}\right) $, where for $t>0$ 
\begin{equation*}
a_{t}\left( p\right) =t^{-m}a\left( tp\right) ,\quad p\in X^{\ast },1\leq
t<\infty .
\end{equation*}%
Let $a_{0}=\varphi a\in \mathcal{C}_{0}^{\infty }\left( X^{\ast }\right) $.
Since 
\begin{equation*}
a\left( p\right) =t^{m}a_{t}\left( p/t\right) ,\quad p\in X^{\ast },1\leq
t<\infty ,
\end{equation*}%
we can write 
\begin{equation*}
a\left( p\right) =a_{0}\left( p\right) +\int_{1}^{\infty }t^{m}\left( \psi
a_{t}\right) \left( \frac{p}{t}\right) \frac{\func{d}t}{t},\quad p\in
X^{\ast },
\end{equation*}%
with the integral also weakly absolutely convergent in $\mathcal{S}^{\ast
}\left( X^{\ast }\right) $. In fact, since $t\leq \left\vert p\right\vert $
in the support of $\psi \left( p/t\right) $ it follows that 
\begin{equation*}
\left\vert t^{m}\left( \psi a_{t}\right) \left( \frac{p}{t}\right)
\right\vert \leq \left( \sup \left\vert \psi a_{t}\right\vert \right)
t^{m-2N}\left\vert p\right\vert ^{2N},\quad p\in X^{\ast },1\leq t<\infty ,
\end{equation*}%
for any $N\in \mathbb{%
%TCIMACRO{\U{2115} }%
%BeginExpansion
\mathbb{N}
%EndExpansion
}$. If we choose $N\in \mathbb{%
%TCIMACRO{\U{2115} }%
%BeginExpansion
\mathbb{N}
%EndExpansion
}$ such that $m<2N$, then we obtain the weakly absolutely convergence in $%
\mathcal{S}^{\ast }\left( X^{\ast }\right) $. Hence 
\begin{equation*}
a=a_{0}+\int_{1}^{\infty }\left( \psi a_{t}\right) _{t^{-1}}\frac{\func{d}t}{%
t}\quad \text{in }\mathcal{S}^{\ast }\left( X^{\ast }\right) .
\end{equation*}%
If we apply the inverse Fourier transformation, $\mathcal{F}^{-1}=\mathcal{F}%
_{X}^{-1}$, to this formula, then we get 
\begin{equation*}
\mathcal{F}^{-1}a=\mathcal{F}^{-1}a_{0}+\int_{1}^{\infty }\mathcal{F}%
^{-1}\left( \psi a_{t}\right) _{t^{-1}}\frac{\func{d}t}{t}
\end{equation*}%
with the integral weakly absolutely convergent in $\mathcal{S}^{\ast }\left(
X\right) $.

We have $\psi a_{t}\in \mathcal{S}\left( X^{\ast }\right) $ and 
\begin{equation*}
\mathcal{F}^{-1}\left( \left( \psi a_{t}\right) _{t^{-1}}\right) \left(
x\right) =t^{m+n}\mathcal{F}^{-1}\left( \psi a_{t}\right) \left( tx\right)
,\quad x\in X.
\end{equation*}

Since the family $\left\{ \mathcal{F}^{-1}\left( \psi a_{t}\right) \right\}
_{1\leq t<\infty }$ is also bounded in $\mathcal{S}\left( X\right) $, it
follows that for any $N\in 
%TCIMACRO{\U{2115} }%
%BeginExpansion
\mathbb{N}
%EndExpansion
$, there is $C_{N}>0$ such that if $x\in X$, then%
\begin{equation*}
\left\vert \mathcal{F}^{-1}\left( \psi a_{t}\right) \left( x\right)
\right\vert \leq C_{N}\left\langle x\right\rangle ^{-N-M},\quad x\in X,1\leq
t<\infty ,
\end{equation*}%
where $M\in 
%TCIMACRO{\U{2115} }%
%BeginExpansion
\mathbb{N}
%EndExpansion
$, $M\geq 1+\max \left\{ 0,m+n\right\} $ is fixed.

It follows that 
\begin{equation}
\left\vert \mathcal{F}^{-1}\left( \left( \psi a_{t}\right) _{t^{-1}}\right)
\left( x\right) \right\vert \leq C_{N}t^{m+n}\left\langle x\right\rangle
^{-N}\left\langle tx\right\rangle ^{-M},\quad x\in X,1\leq t<\infty .
\label{c6}
\end{equation}

We need the following easy consequence of Fubini theorem.

\begin{lemma}
Let $\left( T,\mu \right) $ be a measure space, $\Omega \subset 
%TCIMACRO{\U{211d} }%
%BeginExpansion
\mathbb{R}
%EndExpansion
^{n}$ an open set and $f:\Omega \times T\rightarrow 
%TCIMACRO{\U{2102} }%
%BeginExpansion
\mathbb{C}
%EndExpansion
$ a measurable function.

$(\func{a})$ If for any $\varphi \in \mathcal{C}_{0}^{\infty }\left( \Omega
\right) $ the function $\Omega \times T\ni \left( x,t\right) \rightarrow
\varphi \left( x\right) f\left( x,t\right) \in 
%TCIMACRO{\U{2102} }%
%BeginExpansion
\mathbb{C}
%EndExpansion
$ belongs to $L^{1}\left( \Omega \times T\right) $, then the mapping%
\begin{equation*}
\mathcal{C}_{0}^{\infty }\left( \Omega \right) \ni \varphi \rightarrow \iint
\varphi \left( x\right) f\left( x,t\right) \limfunc{d}x\limfunc{d}\mu \left(
t\right) \in 
%TCIMACRO{\U{2102}}%
%BeginExpansion
\mathbb{C}%
%EndExpansion
\end{equation*}%
define a distribution, the function $\Omega \ni x\rightarrow \int f\left(
x,t\right) \limfunc{d}\mu \left( t\right) \in 
%TCIMACRO{\U{2102} }%
%BeginExpansion
\mathbb{C}
%EndExpansion
$, defined a.e., belongs to $L_{loc}^{1}\left( \Omega \right) $ and we have 
\begin{eqnarray*}
\left\langle \varphi ,\int f\left( \cdot ,t\right) \limfunc{d}\mu \left(
t\right) \right\rangle &=&\iint \varphi \left( x\right) f\left( x,t\right) 
\limfunc{d}x\limfunc{d}\mu \left( t\right) \\
&=&\int \left( \int \varphi \left( x\right) f\left( x,t\right) \limfunc{d}%
x\right) \limfunc{d}\mu \left( t\right) ,\quad \varphi \in \mathcal{C}%
_{0}^{\infty }\left( \Omega \right) .
\end{eqnarray*}

$(\func{b})$ Assume that $\Omega =%
%TCIMACRO{\U{211d} }%
%BeginExpansion
\mathbb{R}
%EndExpansion
^{n}$. If there is $\tau \in 
%TCIMACRO{\U{211d} }%
%BeginExpansion
\mathbb{R}
%EndExpansion
$ such that the function $%
%TCIMACRO{\U{211d} }%
%BeginExpansion
\mathbb{R}
%EndExpansion
^{n}\times T\ni \left( x,t\right) \rightarrow \left\langle x\right\rangle
^{-\tau }f\left( x,t\right) \in 
%TCIMACRO{\U{2102} }%
%BeginExpansion
\mathbb{C}
%EndExpansion
$ belongs to $L^{1}\left( 
%TCIMACRO{\U{211d} }%
%BeginExpansion
\mathbb{R}
%EndExpansion
^{n}\times T\right) $, then the mapping%
\begin{equation*}
\mathcal{S}\left( 
%TCIMACRO{\U{211d} }%
%BeginExpansion
\mathbb{R}
%EndExpansion
^{n}\right) \ni \varphi \rightarrow \iint \varphi \left( x\right) f\left(
x,t\right) \limfunc{d}x\limfunc{d}\mu \left( t\right) \in 
%TCIMACRO{\U{2102}}%
%BeginExpansion
\mathbb{C}%
%EndExpansion
\end{equation*}%
define a temperate distribution, the function $%
%TCIMACRO{\U{211d} }%
%BeginExpansion
\mathbb{R}
%EndExpansion
^{n}\ni x\rightarrow \int f\left( x,t\right) \limfunc{d}\mu \left( t\right)
\in 
%TCIMACRO{\U{2102} }%
%BeginExpansion
\mathbb{C}
%EndExpansion
$, defined a.e., belongs to $L_{loc}^{1}\left( 
%TCIMACRO{\U{211d} }%
%BeginExpansion
\mathbb{R}
%EndExpansion
^{n}\right) $, $\left\langle \cdot \right\rangle ^{-\tau }\int f\left( \cdot
,t\right) \limfunc{d}\mu \left( t\right) \in L^{1}\left( 
%TCIMACRO{\U{211d} }%
%BeginExpansion
\mathbb{R}
%EndExpansion
^{n}\right) $ and we have%
\begin{eqnarray*}
\left\langle \varphi ,\int f\left( \cdot ,t\right) \limfunc{d}\mu \left(
t\right) \right\rangle &=&\iint \varphi \left( x\right) f\left( x,t\right) 
\limfunc{d}x\limfunc{d}\mu \left( t\right) \\
&=&\int \left( \int \varphi \left( x\right) f\left( x,t\right) \limfunc{d}%
x\right) \limfunc{d}\mu \left( t\right) ,\quad \varphi \in \mathcal{S}\left( 
%TCIMACRO{\U{211d} }%
%BeginExpansion
\mathbb{R}
%EndExpansion
^{n}\right) .
\end{eqnarray*}
\end{lemma}

Let us say that a distribution on $X$ is of class $\mathcal{S}$ outside zero
if it is a $\mathcal{C}^{\infty }$ function on $X\backslash \{0\}$ and
decays at infinity, together with all its derivatives, more rapidly than any
power of $\left\langle x\right\rangle ^{-1}$. Then from the representation
formula of $\mathcal{F}^{-1}a$, the estimate $\left( \text{\ref{c6}}\right) $
and part $(\func{a})$ of the previous lemma we conclude that $\mathcal{F}%
^{-1}a$ is of class $\mathcal{S}$ outside zero and we have 
\begin{equation*}
\mathcal{F}^{-1}a\left( x\right) =\mathcal{F}^{-1}a_{0}\left( x\right)
+\int_{1}^{\infty }\mathcal{F}^{-1}\left( \left( \psi a_{t}\right)
_{t^{-1}}\right) \left( x\right) \frac{\func{d}t}{t},\quad x\in X\backslash
\{0\}.
\end{equation*}%
It follows that 
\begin{equation*}
\left\vert \mathcal{F}^{-1}a\left( x\right) \right\vert \leq \left\vert 
\mathcal{F}^{-1}a_{0}\left( x\right) \right\vert +C_{N}\left\langle
x\right\rangle ^{-N}\int_{1}^{\infty }t^{m+n}\left\langle tx\right\rangle
^{-M}\frac{\func{d}t}{t},\quad x\in X\backslash \{0\},
\end{equation*}%
or equivalently 
\begin{equation*}
\left\langle x\right\rangle ^{N}\left\vert \mathcal{F}^{-1}a\left( x\right)
\right\vert \leq \left\langle x\right\rangle ^{N}\left\vert \mathcal{F}%
^{-1}a_{0}\left( x\right) \right\vert +C_{N}\left\vert x\right\vert
^{-m-n}\int_{\left\vert x\right\vert }^{\infty }t^{m+n}\left\langle
t\right\rangle ^{-M}\frac{\func{d}t}{t},\ \ x\in X\backslash \{0\}.
\end{equation*}%
If $m+n>0$, then $C=\int_{0}^{\infty }t^{m+n}\left\langle t\right\rangle
^{-M}\frac{\func{d}t}{t}<\infty $ and%
\begin{equation*}
\left\langle x\right\rangle ^{N}\left\vert \mathcal{F}^{-1}a\left( x\right)
\right\vert \leq \left\langle x\right\rangle ^{N}\left\vert \mathcal{F}%
^{-1}a_{0}\left( x\right) \right\vert +C_{N}C\left\vert x\right\vert
^{-m-n},\quad x\in X\backslash \{0\}.
\end{equation*}

Assume now that $-n<m<0$ and $N\geq n+1$. Then using the estimate $\left( 
\text{\ref{c6}}\right) $ and part $(\func{b})$ of the previous lemma we
conclude that $\mathcal{F}^{-1}a\in L^{1}\left( X\right) $. Since $m_{1}\leq
m_{2}\Rightarrow \mathcal{S}^{m_{1}}\left( X^{\ast }\right) \subset \mathcal{%
S}^{m_{2}}\left( X^{\ast }\right) $, it follows that $\mathcal{F}^{-1}a$
belongs to $L^{1}\left( X\right) $ for any $a\in \bigcup_{m<0}\mathcal{S}%
^{m}\left( X^{\ast }\right) $. Thus we have proved the following

\begin{proposition}
\label{c2}Let $a\in \mathcal{S}^{m}\left( X^{\ast }\right) $. Then:

$(\func{i})$ $\mathcal{F}^{-1}a$ is of class $\mathcal{S}$ outside zero.

$(\func{ii})$ If $m+n>0$, then for any $N\in 
%TCIMACRO{\U{2115} }%
%BeginExpansion
\mathbb{N}
%EndExpansion
$, there is $C_{N}>0$ such that 
\begin{equation*}
\left\vert \mathcal{F}^{-1}a\left( x\right) \right\vert \leq
C_{N}\left\langle x\right\rangle ^{-N}\left( 1+\left\vert x\right\vert
^{-m-n}\right) ,\quad x\in X\backslash \{0\}.
\end{equation*}

$\left( \func{iii}\right) $ If $m<0$, then $\mathcal{F}^{-1}a\in L^{1}\left(
X\right) $.
\end{proposition}

\begin{corollary}
\label{c1}Let $m<-n/2$. If $a\in \mathcal{S}^{m}\left( X^{\ast }\right) $
and $b\in \mathcal{S}^{\infty }\left( X\right) $, then $b\mathcal{F}%
^{-1}a\in L^{2}\left( X\right) $.
\end{corollary}

\begin{proof}
Since $m_{1}\leq m_{2}\Rightarrow \mathcal{S}^{m_{1}}\left( X^{\ast }\right)
\subset \mathcal{S}^{m_{2}}\left( X^{\ast }\right) $, we may suppose that $%
-n<m<-n/2$. Then $\mathcal{F}^{-1}a$ is of class $\mathcal{S}$ outside zero
and for any $N\in 
%TCIMACRO{\U{2115} }%
%BeginExpansion
\mathbb{N}
%EndExpansion
$, there is $C_{N}>0$ such that 
\begin{equation*}
\left\vert \mathcal{F}^{-1}a\left( x\right) \right\vert \leq
C_{N}\left\langle x\right\rangle ^{-N}\left( 1+\left\vert x\right\vert
^{-m-n}\right) ,\quad x\in X\backslash \{0\}.
\end{equation*}%
Since $-m-n>-n/2$ and $N$ can be chosen arbitrarily large, it follows that $b%
\mathcal{F}^{-1}a\in L^{2}\left( X\right) $.
\end{proof}

\begin{corollary}
Let $t>n$ and $m\in \left( n/2,t/2\right) $. If $a\in \mathcal{S}^{-t}\left(
X^{\ast }\right) $ and $b\in \mathcal{S}^{\infty }\left( X\right) $, then $b%
\mathcal{F}^{-1}a\in H^{m}\left( X\right) $.
\end{corollary}

\begin{proof}
Clearly, for any $\varepsilon >0$ there is $M\in 
%TCIMACRO{\U{2115} }%
%BeginExpansion
\mathbb{N}
%EndExpansion
$ such that $\left\langle \cdot \right\rangle ^{-2M}b\in \mathcal{S}%
^{-\varepsilon }\left( X\right) $. It follows that the multiplication
operator by the function by the function $\left\langle \cdot \right\rangle
^{-2M}b=\varphi $ define a bounded operator $M_{\varphi }:H^{r}\left(
X\right) \rightarrow H^{r}\left( X\right) $ for all real $r$. So, it
suffices to show that $\left\langle \cdot \right\rangle ^{2M}\mathcal{F}%
^{-1}a\in H^{m}\left( X\right) $ for any $M\in 
%TCIMACRO{\U{2115} }%
%BeginExpansion
\mathbb{N}
%EndExpansion
$. Since $\left( 1-\triangle \right) ^{m/2}\left[ \left\langle \cdot
\right\rangle ^{2M}\mathcal{F}^{-1}a\right] =\mathcal{F}^{-1}c$ with $%
c=\left\langle \cdot \right\rangle ^{m}\left( 1-\triangle _{X^{\ast
}}\right) ^{M}a\in \mathcal{S}^{m-t}\left( X^{\ast }\right) $ and $%
m-t<-t/2<-n/2,$ the previous corollary implies that $\left( 1-\triangle
\right) ^{m/2}\left[ \left\langle \cdot \right\rangle ^{2M}\mathcal{F}^{-1}a%
\right] \in L^{2}\left( X\right) $. Hence $\left\langle \cdot \right\rangle
^{2M}\mathcal{F}^{-1}a\in H^{m}\left( X\right) $.
\end{proof}

\begin{remark}
Let $V$ be an euclidean space, $A\in \limfunc{End}_{%
%TCIMACRO{\U{211d} }%
%BeginExpansion
\mathbb{R}
%EndExpansion
}\left( V\right) $ and $\chi \in \mathcal{S}\left( V\right) $.

$(\func{i})$ If $A$ is invvertible, then $\left\langle v\right\rangle \leq
\max \left\{ 1,\left\Vert A^{-1}\right\Vert \right\} \left\langle
Av\right\rangle $, $v\in V.$

$(\func{ii})$ If $h\in 
%TCIMACRO{\U{211d} }%
%BeginExpansion
\mathbb{R}
%EndExpansion
$ is such that $\left\vert h\right\vert \left\Vert A\right\Vert \leq 1/2$,
then $\left\langle v\right\rangle \leq 2\left\langle v+\lambda
hAv\right\rangle $, $v\in V$, $0\leq \lambda \leq 1$ and by Taylor's formula%
\begin{equation*}
\left\langle v\right\rangle ^{k}\left\vert \chi \left( v+hAv\right) -\chi
\left( v\right) \right\vert \leq 2^{k+1}\left\Vert A\right\Vert \left\vert
h\right\vert \int_{0}^{1}\left\langle v+\lambda hAv\right\rangle
^{k+1}\left\Vert \chi ^{\prime }\left( v+\lambda hAv\right) \right\Vert 
\limfunc{d}\lambda ,
\end{equation*}%
for $v\in V$.
\end{remark}

Let $\tau \in \limfunc{End}_{%
%TCIMACRO{\U{211d} }%
%BeginExpansion
\mathbb{R}
%EndExpansion
}\left( X\right) $. We consider $C_{\tau }\in \limfunc{End}_{%
%TCIMACRO{\U{211d} }%
%BeginExpansion
\mathbb{R}
%EndExpansion
}\left( X\times X\right) $ defined by $C_{\tau }\left( x,y\right) =\left(
\left( 1-\tau \right) x+\tau y,x-y\right) $, $\left( x,y\right) \in X\times
X $ with the inverse $C_{\tau }^{-1}$ given by $C_{\tau }^{-1}\left(
v,u\right) =\left( v+\tau u,v-\left( 1-\tau \right) u\right) $, $\left(
u,v\right) \in X\times X$. If $\tau ,\tau _{0}\in \limfunc{End}_{%
%TCIMACRO{\U{211d} }%
%BeginExpansion
\mathbb{R}
%EndExpansion
}\left( X\right) $, then $C_{\tau }\circ C_{\tau _{0}}^{-1}=1_{X\times
X}+\left( \tau _{0}-\tau \right) A$, where $A\left( v,u\right) =\left(
u,0\right) $ satisfies $\left\Vert A\right\Vert =1$.

Let $\varphi \in \mathcal{S}\left( X\times X\right) $. Then by the above
remark it follows that for any $k\in 
%TCIMACRO{\U{2115} }%
%BeginExpansion
\mathbb{N}
%EndExpansion
$ there is $C_{\varphi ,k}>0$ such that 
\begin{equation*}
\left\langle \left( x,y\right) \right\rangle ^{k}\left\vert \left( \varphi
\circ C_{\tau }-\varphi \circ C_{\tau _{0}}\right) \left( x,y\right)
\right\vert \leq C_{\varphi ,k}\left\vert \tau -\tau _{0}\right\vert
\end{equation*}%
if $\left\vert \tau -\tau _{0}\right\vert \leq 1/2$.

Using this estimate we obtain a simple but useful lemma.

\begin{lemma}
\label{c4}Let $a\in \mathcal{S}\left( X^{\ast }\right) $, $b\in \mathcal{S}%
\left( X\right) $ and $c\in \mathcal{S}^{\infty }\left( X\times X\right) $.
For $\tau \in \limfunc{End}_{%
%TCIMACRO{\U{211d} }%
%BeginExpansion
\mathbb{R}
%EndExpansion
}\left( X\right) $ we define $\mathcal{K}\left( \cdot ,\cdot ;\tau \right)
=c\cdot \left[ \left( \mathcal{F}^{-1}a\otimes b\right) \circ C_{\tau }%
\right] $ i.e.%
\begin{equation*}
\mathcal{K}\left( x,y;\tau \right) =c\left( x,y\right) b\left( x-y\right) 
\mathcal{F}^{-1}a\left( \left( 1-\tau \right) x+\tau y\right) ,\quad \left(
x,y\right) \in X\times X.
\end{equation*}%
Then $\mathcal{K}\left( \cdot ,\cdot ;\tau \right) \in L^{2}\left( X\times
X\right) $ and the mapping%
\begin{equation*}
\limfunc{End}\nolimits_{%
%TCIMACRO{\U{211d} }%
%BeginExpansion
\mathbb{R}
%EndExpansion
}\left( X\right) \ni \tau \rightarrow \mathcal{K}\left( \cdot ,\cdot ;\tau
\right) \in L^{2}\left( X\times X\right)
\end{equation*}%
is continuous.
\end{lemma}

Let $\mathcal{U}_{0}$ be the set of all invetible endomorphisms of $X$, $%
\mathcal{U}_{1}=1_{X}+\mathcal{U}_{0}$ and $\mathcal{U}=\mathcal{U}_{0}\cup 
\mathcal{U}_{1}$. It is clear that all these sets are open in $\limfunc{End}%
_{%
%TCIMACRO{\U{211d} }%
%BeginExpansion
\mathbb{R}
%EndExpansion
}\left( X\right) $ and $%
%TCIMACRO{\U{211d} }%
%BeginExpansion
\mathbb{R}
%EndExpansion
\cdot 1_{X}\subset \mathcal{U}$. We shall extend the results in \cite{Arsu}
from the case when $\tau $ is a real number to the case when $\tau $ belongs
to the open subset $\mathcal{U}\subset \limfunc{End}_{%
%TCIMACRO{\U{211d} }%
%BeginExpansion
\mathbb{R}
%EndExpansion
}\left( X\right) $.

It is useful to estimate the norm of of an element $f\in H^{m}\left(
X\right) $ without referring to the Fourier transform of $f$. If $m=\left[ m%
\right] +\mu $ and $r>0$, then 
\begin{eqnarray*}
\left\Vert f\right\Vert _{H^{m}\left( X\right) }^{2} &\approx
&\sum_{\left\vert \alpha \right\vert \leq \left[ m\right] }\left\Vert
\partial ^{\alpha }f\right\Vert _{L^{2}\left( X\right)
}^{2}+\sum_{\left\vert \alpha \right\vert =\left[ m\right]
}\iint_{\left\vert z\right\vert \leq r}\frac{\left\vert \partial ^{\alpha
}f\left( x+z\right) -\partial ^{\alpha }f\left( x\right) \right\vert ^{2}}{%
\left\vert z\right\vert ^{n+2\mu }}\limfunc{d}x\limfunc{d}z \\
&\equiv &\left\Vert f\right\Vert _{m,r}^{2}.
\end{eqnarray*}

\begin{lemma}
Let $s,t>n$ and $m\in \left( n/2,t/2\right) $. Let $a\in \mathcal{S}%
^{-t}\left( X^{\ast }\right) $, $b\in \mathcal{S}^{-s}\left( X\right) $ and $%
c\in \mathcal{S}^{s/2}\left( X\right) $. For $\tau \in \limfunc{End}_{%
%TCIMACRO{\U{211d} }%
%BeginExpansion
\mathbb{R}
%EndExpansion
}\left( X\right) $ we put $\mathcal{K}_{a,b}\left( \cdot ,\cdot ;\tau
\right) =\left( \mathcal{F}^{-1}a\otimes b\right) \circ C_{\tau }$ i.e. 
\begin{equation*}
\mathcal{K}_{a,b}\left( x,y;\tau \right) =b\left( x-y\right) \mathcal{F}%
^{-1}a\left( \left( 1-\tau \right) x+\tau y\right) ,\quad \left( x,y\right)
\in X\times X.
\end{equation*}

$(\func{i})$ If $\tau \in \mathcal{U}_{1}$, then the function%
\begin{equation*}
X\times X\ni \left( x,y\right) \rightarrow \left( 1-\triangle \right)
^{m/2}\left( c\left( \cdot \right) \mathcal{K}_{a,b}\left( \cdot ,y;\tau
\right) \right) \left( x\right) \equiv \mathcal{K}_{a,b,c}^{1}\left(
x,y;\tau \right) \in 
%TCIMACRO{\U{2102}}%
%BeginExpansion
\mathbb{C}%
%EndExpansion
\end{equation*}%
is in $L^{2}\left( X\times X\right) $ and the mapping%
\begin{equation*}
\mathcal{U}_{1}\ni \tau \rightarrow \mathcal{K}_{a,b,c}^{1}\left( \cdot
,\cdot ;\tau \right) \in L^{2}\left( X\times X\right)
\end{equation*}%
is continuous.

$(\func{ii})$ If $\tau \in \mathcal{U}_{0}$, then the function 
\begin{equation*}
X\times X\ni \left( x,y\right) \rightarrow \left( 1-\triangle \right)
^{m/2}\left( c\left( \cdot \right) \mathcal{K}_{a,b}\left( x,\cdot ;\tau
\right) \right) \left( y\right) \equiv \mathcal{K}_{a,b,c}^{0}\left(
x,y;\tau \right) \in 
%TCIMACRO{\U{2102}}%
%BeginExpansion
\mathbb{C}%
%EndExpansion
\end{equation*}%
is in $L^{2}\left( X\times X\right) $ and the mapping%
\begin{equation*}
\mathcal{U}_{0}\ni \tau \rightarrow \mathcal{K}_{a,b,c}^{0}\left( \cdot
,\cdot ;\tau \right) \in L\left( X\times X\right)
\end{equation*}%
is continuous.
\end{lemma}

\begin{proof}
Let $\psi =\mathcal{F}^{-1}a\in H^{m}\left( X\right) \cap L^{1}\left(
X\right) $. We observe that 
\begin{multline*}
\left\Vert \mathcal{K}_{a,b,c}^{1}\left( \cdot ,\cdot ;\tau \right)
\right\Vert _{L^{2}}^{2} \\
=\iint \left\vert \left( 1-\triangle \right) ^{m/2}\left( c\left( \cdot
\right) b\left( \cdot -y\right) \psi \left( \left( 1-\tau \right) \cdot
+\tau y\right) \right) \left( x\right) \right\vert ^{2}\limfunc{d}x\limfunc{d%
}y \\
=\int \left\Vert c\left( \cdot \right) b\left( \cdot -y\right) \psi \left(
\left( 1-\tau \right) \cdot +\tau y\right) \right\Vert _{H^{m}\left(
X\right) }^{2}\limfunc{d}y \\
\leq C\sum_{\left\vert \alpha \right\vert \leq \left[ m\right] }\int
\left\Vert f_{\alpha }\left( \cdot ,y\right) \right\Vert _{L^{2}\left(
X\right) }^{2}\limfunc{d}y \\
+C\sum_{\left\vert \alpha \right\vert =\left[ m\right] }\iiint_{\left\vert
z\right\vert \leq 1}\frac{\left\vert f_{\alpha }\left( x+z,y\right)
-f_{\alpha }\left( x,y\right) \right\vert ^{2}}{\left\vert z\right\vert
^{n+2\mu }}\limfunc{d}x\limfunc{d}y\limfunc{d}z,
\end{multline*}%
where we set $f_{\alpha }\left( x,y\right) =\partial ^{\alpha }\left(
c\left( \cdot \right) b\left( \cdot -y\right) \psi \left( \left( 1-\tau
\right) \cdot +\tau y\right) \right) \left( x\right) $. When $m\in 
%TCIMACRO{\U{2115} }%
%BeginExpansion
\mathbb{N}
%EndExpansion
$ then the sum $\sum_{\left\vert \alpha \right\vert =\left[ m\right]
}\iiint_{\left\vert z\right\vert \leq 1}...$ don't appear in the above
estimate.

For $1\leq \left\vert \alpha \right\vert \leq \left[ m\right] $, Leibniz'
formula implies the equality%
\begin{equation*}
f_{\alpha }\left( x,y\right) =\sum_{\beta +\gamma +\delta =\alpha }\frac{%
\alpha !}{\beta !\gamma !\delta !}P_{\delta }\left( \tau \right) \partial
^{\beta }c\left( x\right) \partial ^{\gamma }b\left( x-y\right) \partial
^{\delta }\psi \left( \left( 1-\tau \right) x+\tau y\right)
\end{equation*}%
where $P_{\delta }\left( \tau \right) $ is a polynomial of degree $%
\left\vert \delta \right\vert $ in $\tau _{ij}$.

Let $K$ be a compact subset of $\mathcal{U}_{1}$. We have to estimate
several types of terms.

We begin with the simplest type. For $\tau \in K$ we have 
\begin{multline*}
\iint \left\vert \partial ^{\beta }c\left( x\right) \partial ^{\gamma
}b\left( x-y\right) \partial ^{\delta }\psi \left( \left( 1-\tau \right)
x+\tau y\right) \right\vert ^{2}\limfunc{d}x\limfunc{d}y \\
=\iint \left\vert \partial ^{\beta }c\left( v+\tau u\right) \partial
^{\gamma }b\left( u\right) \partial ^{\delta }\psi \left( v\right)
\right\vert ^{2}\limfunc{d}u\limfunc{d}v \\
\leq C\left( \beta ,\gamma ,\delta \right) \iint \left\langle v+\tau
u\right\rangle ^{s}\left\langle u\right\rangle ^{-2s}\left\vert \partial
^{\delta }\psi \left( v\right) \right\vert ^{2}\limfunc{d}u\limfunc{d}v \\
\leq C\left( K,\beta ,\gamma ,\delta ,s\right) \iint \left\langle
u\right\rangle ^{-s}\left\langle v\right\rangle ^{s}\left\vert \partial
^{\delta }\psi \left( v\right) \right\vert ^{2}\limfunc{d}u\limfunc{d}v \\
=C\left( K,\beta ,\gamma ,\delta ,s\right) \left( \int \left\langle
u\right\rangle ^{-s}\limfunc{d}u\right) \left( \int \left\langle
v\right\rangle ^{s}\left\vert \partial ^{\delta }\psi \left( v\right)
\right\vert ^{2}\limfunc{d}v\right) <\infty .
\end{multline*}%
Here we used the change of variable $\left\{ 
\begin{array}{l}
u=x-y \\ 
v=\left( 1-\tau \right) x+\tau y%
\end{array}%
\right. $, Peetre's inequality $\left\langle x+y\right\rangle ^{s}\leq
2^{\left\vert s\right\vert /2}\left\langle x\right\rangle ^{\left\vert
s\right\vert }\left\langle y\right\rangle ^{s}$ and Corollary \ref{c1}. $%
\partial ^{\delta }\psi =i^{\left\vert \delta \right\vert }\mathcal{F}%
^{-1}\left( p^{\delta }a\right) $ and $p^{\delta }a\in \mathcal{S}%
^{\left\vert \delta \right\vert -t}\left( X^{\ast }\right) $ with $%
\left\vert \delta \right\vert -t\leq m-t<-n/2$.

For $\left\vert \alpha \right\vert =\left[ m\right] $ and $\alpha =\beta
+\gamma +\delta $ we set 
\begin{eqnarray*}
f\left( x,y\right) &=&f_{\beta ,\gamma ,\delta }\left( x,y\right) =\partial
^{\beta }c\left( x\right) \partial ^{\gamma }b\left( x-y\right) \partial
^{\delta }\psi \left( \left( 1-\tau \right) x+\tau y\right) \\
&=&h\left( x-y,\left( 1-\tau \right) x+\tau y\right) \partial ^{\delta }\psi
\left( \left( 1-\tau \right) x+\tau y\right) ,
\end{eqnarray*}%
where $h\left( u,v\right) =\partial ^{\beta }c\left( v+\tau u\right)
\partial ^{\gamma }b\left( u\right) $.

The second type of terms is given by the integral 
\begin{multline*}
I=\iint\limits_{X\times X}\int\limits_{\left\vert z\right\vert \leq 1}\frac{%
\left\vert f\left( x+z,y\right) -f\left( x,y\right) \right\vert ^{2}}{%
\left\vert z\right\vert ^{n+2\mu }}\limfunc{d}x\limfunc{d}y\limfunc{d}z \\
=\int\limits_{\left\vert z\right\vert \leq 1}\iint\limits_{X\times X}\frac{%
\left\vert f\left( v+\tau u+z,v-\left( 1-\tau \right) u\right) -f\left(
v+\tau u,v-\left( 1-\tau \right) u\right) \right\vert ^{2}}{\left\vert
z\right\vert ^{n+2\mu }}\limfunc{d}u\limfunc{d}v\limfunc{d}z \\
=\int\limits_{\left\vert z\right\vert \leq 1}\iint\limits_{X\times X}\frac{%
\left\vert h\left( u+z,v+\left( 1-\tau \right) z\right) \partial ^{\delta
}\psi \left( v+\left( 1-\tau \right) z\right) -h\left( u,v\right) \partial
^{\delta }\psi \left( v\right) \right\vert ^{2}}{\left\vert z\right\vert
^{n+2\mu }}\limfunc{d}u\limfunc{d}v\limfunc{d}z.
\end{multline*}%
Here we used again the change of variable $\left\{ 
\begin{array}{l}
u=x-y \\ 
v=\left( 1-\tau \right) x+\tau y%
\end{array}%
\right. $. We have 
\begin{multline*}
\int\limits_{\left\vert z\right\vert \leq 1}\iint\limits_{X\times X}\frac{%
\left\vert h\left( u+z,v+\left( 1-\tau \right) z\right) \partial ^{\delta
}\psi \left( v+\left( 1-\tau \right) z\right) -h\left( u,v\right) \partial
^{\delta }\psi \left( v\right) \right\vert ^{2}}{\left\vert z\right\vert
^{n+2\mu }}\limfunc{d}u\limfunc{d}v\limfunc{d}z\quad \\
\leq I_{1}+I_{2},
\end{multline*}%
where 
\begin{equation*}
I_{1}=\int\limits_{\left\vert z\right\vert \leq 1}\iint\limits_{X\times X}%
\frac{\left\vert h\left( u+z,v+\left( 1-\tau \right) z\right) -h\left(
u,v\right) \right\vert ^{2}\left\vert \partial ^{\delta }\psi \left(
v\right) \right\vert ^{2}}{\left\vert z\right\vert ^{n+2\mu }}\limfunc{d}u%
\limfunc{d}v\limfunc{d}z,
\end{equation*}%
and%
\begin{equation*}
I_{2}=\int\limits_{\left\vert z\right\vert \leq 1}\iint\limits_{X\times X}%
\frac{\left\vert h\left( u+z,v+\left( 1-\tau \right) z\right) \right\vert
^{2}\left\vert \partial ^{\delta }\psi \left( v+\left( 1-\tau \right)
z\right) -\partial ^{\delta }\psi \left( v\right) \right\vert ^{2}}{%
\left\vert z\right\vert ^{n+2\mu }}\limfunc{d}u\limfunc{d}v\limfunc{d}z.
\end{equation*}

The estimate of $I_{1}$ is easier. First we observe that 
\begin{eqnarray*}
h\left( u+z,v+\left( 1-\tau \right) z\right) &=&\partial ^{\beta }c\left(
v+\left( 1-\tau \right) z+\tau u+\tau z\right) \partial ^{\gamma }b\left(
u+z\right) \\
&=&\partial ^{\beta }c\left( v+\tau u+z\right) \partial ^{\gamma }b\left(
u+z\right) .
\end{eqnarray*}%
Then we use the mean value theorem and Peetre's inequality to estimate the
difference. For $\tau \in K$ we have%
\begin{eqnarray*}
\left\vert h\left( u+z,v+\left( 1-\tau \right) z\right) -h\left( u,v\right)
\right\vert ^{2} &\leq &C\left\vert z\right\vert ^{2}\sup_{0\leq \lambda
\leq 1}\left\langle v+\tau u+\lambda z\right\rangle ^{s}\left\langle
u+\lambda z\right\rangle ^{-2s} \\
&\leq &C\left( K,s\right) \left\vert z\right\vert ^{2}\left\langle
v\right\rangle ^{s}\left\langle u\right\rangle ^{-s}.
\end{eqnarray*}%
Hence 
\begin{equation*}
I_{1}\leq C\left( K,s\right) \left( \int_{X}\left\langle u\right\rangle ^{-s}%
\limfunc{d}u\right) \left( \int_{X}\left\langle v\right\rangle
^{s}\left\vert \partial ^{\delta }\psi \left( v\right) \right\vert ^{2}%
\limfunc{d}v\right) \left( \int_{\left\vert z\right\vert \leq 1}\frac{%
\limfunc{d}z}{\left\vert z\right\vert ^{n+2\mu -2}}\right) <\infty
\end{equation*}%
since Corollary \ref{c1} is applicable. $\mathcal{F}\left( \partial ^{\delta
}\psi \right) =i^{\left\vert \delta \right\vert }p^{\delta }a$ and $%
p^{\delta }a\in \mathcal{S}^{\left\vert \delta \right\vert -t}\left( X^{\ast
}\right) $ with $\left\vert \delta \right\vert -t\leq m-t<-n/2$.

Let now estimate $I_{2}$. By using again Peetre's inequality, the fact that $%
\tau $ belongs to $K$, a compact subset of $\mathcal{U}_{1}$, and the fact
that $\left\vert z\right\vert \leq 1$, it follows that%
\begin{eqnarray*}
\left\vert h\left( u+z,v+\left( 1-\tau \right) z\right) \right\vert ^{2}
&=&\left\vert \partial ^{\beta }c\left( v+\tau u+z\right) \partial ^{\gamma
}b\left( u+z\right) \right\vert ^{2} \\
&\leq &C\left\langle v+\tau u+z\right\rangle ^{s}\left\langle
u+z\right\rangle ^{-2s}\leq C\left( K,s\right) \left\langle v\right\rangle
^{s}\left\langle u\right\rangle ^{-s}.
\end{eqnarray*}%
Hence 
\begin{multline*}
I_{2}=\int\limits_{\left\vert z\right\vert \leq 1}\iint\limits_{X\times X}%
\frac{\left\vert h\left( u+z,v+\left( 1-\tau \right) z\right) \right\vert
^{2}\left\vert \partial ^{\delta }\psi \left( v+\left( 1-\tau \right)
z\right) -\partial ^{\delta }\psi \left( v\right) \right\vert ^{2}}{%
\left\vert z\right\vert ^{n+2\mu }}\limfunc{d}u\limfunc{d}v\limfunc{d}z \\
\leq C\left( K,s\right) \left( \int_{X}\left\langle u\right\rangle ^{-s}%
\limfunc{d}u\right) \left( \int_{X}\int_{\left\vert z\right\vert \leq 1}%
\frac{\left\langle v\right\rangle ^{s}\left\vert \partial ^{\delta }\psi
\left( v+\left( 1-\tau \right) z\right) -\partial ^{\delta }\psi \left(
v\right) \right\vert ^{2}}{\left\vert z\right\vert ^{n+2\mu }}\limfunc{d}v%
\limfunc{d}z\right)
\end{multline*}%
so it remains to evaluate the integral 
\begin{equation*}
J=\int_{X}\int_{\left\vert z\right\vert \leq 1}\frac{\left\langle
v\right\rangle ^{s}\left\vert \partial ^{\delta }\psi \left( v+\left( 1-\tau
\right) z\right) -\partial ^{\delta }\psi \left( v\right) \right\vert ^{2}}{%
\left\vert z\right\vert ^{n+2\mu }}\limfunc{d}v\limfunc{d}z.
\end{equation*}%
Now we shall use the fact that $1-\tau $ is invertible if $\tau \in \mathcal{%
U}_{1}$. We have 
\begin{eqnarray*}
J &\leq &\left\Vert 1-\tau \right\Vert ^{n+2\mu }\int_{X}\int_{\left\vert
z\right\vert \leq 1}\frac{\left\langle v\right\rangle ^{s}\left\vert
\partial ^{\delta }\psi \left( v+\left( 1-\tau \right) z\right) -\partial
^{\delta }\psi \left( v\right) \right\vert ^{2}}{\left\vert \left( 1-\tau
\right) z\right\vert ^{n+2\mu }}\limfunc{d}v\limfunc{d}z \\
&\leq &\frac{\left\Vert 1-\tau \right\Vert ^{n+2\mu }}{\left\vert \det
\left( 1-\tau \right) \right\vert }\int_{X}\int_{\left\vert \zeta
\right\vert \leq \left\Vert 1-\tau \right\Vert }\frac{\left\langle
v\right\rangle ^{s}\left\vert \partial ^{\delta }\psi \left( v+\zeta \right)
-\partial ^{\delta }\psi \left( v\right) \right\vert ^{2}}{\left\vert \zeta
\right\vert ^{n+2\mu }}\limfunc{d}v\limfunc{d}\zeta
\end{eqnarray*}%
Next we shall split $X$ in two regions $X=\left\{ v:\left\vert v\right\vert
\geq 2\left\Vert 1-\tau \right\Vert \right\} \cup \left\{ v:\left\vert
v\right\vert \leq 2\left\Vert 1-\tau \right\Vert \right\} $. Then 
\begin{equation*}
\int_{X}\int_{\left\vert \zeta \right\vert \leq \left\Vert 1-\tau
\right\Vert }\frac{\left\langle v\right\rangle ^{s}\left\vert \partial
^{\delta }\psi \left( v+\zeta \right) -\partial ^{\delta }\psi \left(
v\right) \right\vert ^{2}}{\left\vert \zeta \right\vert ^{n+2\mu }}\limfunc{d%
}v\limfunc{d}\zeta =J_{1}+J_{2}
\end{equation*}%
where 
\begin{eqnarray*}
J_{1} &=&\int_{\left\vert v\right\vert \geq 2\left\Vert 1-\tau \right\Vert
}\int_{\left\vert \zeta \right\vert \leq \left\Vert 1-\tau \right\Vert }%
\frac{\left\langle v\right\rangle ^{s}\left\vert \partial ^{\delta }\psi
\left( v+\zeta \right) -\partial ^{\delta }\psi \left( v\right) \right\vert
^{2}}{\left\vert \zeta \right\vert ^{n+2\mu }}\limfunc{d}v\limfunc{d}\zeta ,
\\
J_{2} &=&\int_{\left\vert v\right\vert \leq 2\left\Vert 1-\tau \right\Vert
}\int_{\left\vert \zeta \right\vert \leq \left\Vert 1-\tau \right\Vert }%
\frac{\left\langle v\right\rangle ^{s}\left\vert \partial ^{\delta }\psi
\left( v+\zeta \right) -\partial ^{\delta }\psi \left( v\right) \right\vert
^{2}}{\left\vert \zeta \right\vert ^{n+2\mu }}\limfunc{d}v\limfunc{d}\zeta .
\end{eqnarray*}

The estimate of $J_{1}$. If $\left\vert v\right\vert \geq 2\left\Vert 1-\tau
\right\Vert $, $\left\vert \zeta \right\vert \leq \left\Vert 1-\tau
\right\Vert $ and $0\leq \lambda \leq 1$, then $\left\vert v+\lambda \zeta
\right\vert \geq \left\vert v\right\vert -\left\vert \lambda \zeta
\right\vert \geq \left\Vert 1-\tau \right\Vert $. In view of Proposition \ref%
{c2} $(\func{i})$, it follows that for any $N\in 
%TCIMACRO{\U{2115} }%
%BeginExpansion
\mathbb{N}
%EndExpansion
$ there is $C_{N}>0$ such that 
\begin{eqnarray*}
\left\vert \partial ^{\delta }\psi \left( v+\zeta \right) -\partial ^{\delta
}\psi \left( v\right) \right\vert ^{2} &=&\left\vert \left\langle \zeta
,\int_{0}^{1}\left( \partial ^{\delta }\psi \right) ^{\prime }\left(
v+\lambda \zeta \right) \limfunc{d}\lambda \right\rangle \right\vert ^{2} \\
&\leq &C_{N}\left\vert \zeta \right\vert ^{2}\sup_{0\leq \lambda \leq
1}\left\langle v+\lambda \zeta \right\rangle ^{-N} \\
&\leq &C_{N}\left( K\right) \left\vert \zeta \right\vert ^{2}\left\langle
v\right\rangle ^{-N}.
\end{eqnarray*}%
If we take $N=\left[ s\right] +1+n+1$ we obtain that 
\begin{eqnarray*}
J_{1} &\leq &C_{N}\left( K\right) \left\Vert 1-\tau \right\Vert ^{2\left(
1-\mu \right) }\left( \int_{X}\left\langle v\right\rangle ^{-n-1}\limfunc{d}%
v\right) \left( \int_{\left\vert \zeta \right\vert \leq 1}\frac{\limfunc{d}%
\zeta }{\left\vert \zeta \right\vert ^{n+2\mu -2}}\right) \\
&\leq &C\left( K,s\right) \left( \int_{X}\left\langle v\right\rangle ^{-n-1}%
\limfunc{d}v\right) \left( \int_{\left\vert \zeta \right\vert \leq 1}\frac{%
\limfunc{d}\zeta }{\left\vert \zeta \right\vert ^{n+2\mu -2}}\right) <\infty
.
\end{eqnarray*}

The estimate of $J_{2}$. In order to estimate $J_{2}$ we shall use the fact
that $\psi =\mathcal{F}^{-1}a\in H^{m}\left( X\right) $. For $\tau \in K$ we
have%
\begin{eqnarray*}
J_{2} &=&\int_{\left\vert v\right\vert \leq 2\left\Vert 1-\tau \right\Vert
}\int_{\left\vert \zeta \right\vert \leq \left\Vert 1-\tau \right\Vert }%
\frac{\left\langle v\right\rangle ^{s}\left\vert \partial ^{\delta }\psi
\left( v+\zeta \right) -\partial ^{\delta }\psi \left( v\right) \right\vert
^{2}}{\left\vert \zeta \right\vert ^{n+2\mu }}\limfunc{d}v\limfunc{d}\zeta \\
&\leq &\left( 1+4\left\Vert 1-\tau \right\Vert ^{2}\right)
^{s/2}\int_{X}\int_{\left\vert \zeta \right\vert \leq \left\Vert 1-\tau
\right\Vert }\frac{\left\vert \partial ^{\delta }\psi \left( v+\zeta \right)
-\partial ^{\delta }\psi \left( v\right) \right\vert ^{2}}{\left\vert \zeta
\right\vert ^{n+2\mu }}\limfunc{d}v\limfunc{d}\zeta \\
&\leq &C\left( K\right) \left\Vert \psi \right\Vert _{H^{\left\vert \delta
\right\vert +\mu }\left( X\right) }^{2}\leq C\left( K\right) \left\Vert \psi
\right\Vert _{H^{m}\left( X\right) }^{2}<\infty .
\end{eqnarray*}

Now pass to the proof of the continuity of the mapping 
\begin{equation*}
\mathcal{U}_{1}\ni \tau \rightarrow \mathcal{K}_{a,b,c}^{1}\left( \cdot
,\cdot ;\tau \right) \in L^{2}\left( X\times X\right) .
\end{equation*}

Let $K$ be a compact subset of $\mathcal{U}_{1}$. We have shown that the
restriction of the above mapping to $K$ is in $L^{\infty }\left(
K;L^{2}\left( X\times X\right) \right) $.

Now we observe that the bilinear mapping%
\begin{equation*}
\mathcal{S}^{-t}\left( X^{\ast }\right) \times \mathcal{S}^{-s}\left(
X\right) \ni \left( a,b\right) \rightarrow \mathcal{K}_{a,b,c}^{1}\in
L^{\infty }\left( K;L^{2}\left( X\times X\right) \right)
\end{equation*}%
is separately continuous. We shall prove only the continuity in the first
variable since the continuity in the second variable can be done in the same
manner. To show that the bilinear mapping $\left( a,b\right) \rightarrow 
\mathcal{K}_{a,b,c}^{1}$ is continuous in the first variable, we shall use
the closed graph theorem. If $a_{n}\rightarrow a$ in $\mathcal{S}^{-t}\left(
X^{\ast }\right) $ and $\mathcal{K}_{a_{n},b,c}^{1}\rightarrow \mathcal{K}$
in $L^{\infty }\left( K;L^{2}\left( X\times X\right) \right) $, then $%
\mathcal{K}_{a_{n},b,c}^{1}\left( \tau \right) \rightarrow \mathcal{K}%
_{a,b,c}^{1}\left( \tau \right) $ in $\mathcal{S}^{\ast }\left( X\times
X\right) $ for every $\tau \in K$ and $\mathcal{K}_{a_{n},b,c}^{1}\left(
\tau \right) \rightarrow \mathcal{K}\left( \tau \right) $ in $L^{2}\left(
X\times X\right) $ a.e. $\tau \in K$. It follows that $\mathcal{K}%
_{a,b,c}^{1}=\mathcal{K}$ in $L^{\infty }\left( K;L^{2}\left( X\times
X\right) \right) $.

But a bilinear mapping in the product of a Fr\'{e}chet space and a
metrizable space is continuous if it is separately continuous. Hence the
bilinear mapping $\left( a,b\right) \rightarrow \mathcal{K}_{a,b,c}^{1}$ is
continuous.

Let $t^{\prime }\in \left( 2m,t\right) $ and $s^{\prime }\in \left(
n,s\right) $. Then there are two sequences $\left\{ a_{n}\right\} \subset 
\mathcal{S}\left( X^{\ast }\right) $ and $\left\{ b_{n}\right\} \subset 
\mathcal{S}\left( X\right) $ such that $a_{n}\rightarrow a$ in $\mathcal{S}%
^{-t^{\prime }}\left( X^{\ast }\right) $ and $b_{n}\rightarrow b$ in $%
\mathcal{S}^{-s^{\prime }}\left( X\right) $ in virtue of Corollary \ref{c3}.
Since the mapping%
\begin{equation*}
\mathcal{S}^{-t^{\prime }}\left( X^{\ast }\right) \times \mathcal{S}%
^{-s^{\prime }}\left( X\right) \ni \left( a,b\right) \rightarrow \mathcal{K}%
_{a,b,c}^{1}\in L^{\infty }\left( K;L^{2}\left( X\times X\right) \right)
\end{equation*}%
is continuous, $\lim \left\Vert \mathcal{K}_{a_{n},b_{n},c}^{1}-\mathcal{K}%
_{a,b,c}^{1}\right\Vert _{L^{\infty }}=0$. Hence we have to prove the
continuity only when $a\in \mathcal{S}\left( X^{\ast }\right) $ and $b\in 
\mathcal{S}\left( X\right) $.

If $m\leq M$, $M\in 
%TCIMACRO{\U{2115} }%
%BeginExpansion
\mathbb{N}
%EndExpansion
$, then by Leibniz' rule we obtain 
\begin{multline*}
\left\Vert \mathcal{K}_{a,b,c}^{1}\left( \cdot ,\cdot ;\tau \right) -%
\mathcal{K}_{a,b,c}^{1}\left( \cdot ,\cdot ;\tau _{0}\right) \right\Vert
_{L^{2}\left( X\times X\right) }^{2} \\
=\int \left\Vert c\left( \cdot \right) b\left( \cdot -y\right) \left( 
\mathcal{F}^{-1}a\left( \left( 1-\tau \right) \cdot +\tau y\right) -\mathcal{%
F}^{-1}a\left( \left( 1-\tau _{0}\right) \cdot +\tau _{0}y\right) \right)
\right\Vert _{H^{m}\left( X\right) }^{2}\limfunc{d}y \\
\leq \int \left\Vert c\left( \cdot \right) b\left( \cdot -y\right) \left( 
\mathcal{F}^{-1}a\left( \left( 1-\tau \right) \cdot +\tau y\right) -\mathcal{%
F}^{-1}a\left( \left( 1-\tau _{0}\right) \cdot +\tau _{0}y\right) \right)
\right\Vert _{H^{M}\left( X\right) }^{2}\limfunc{d}y \\
\leq C\left( K,M\right) \sum_{\left\vert \alpha \right\vert +\left\vert
\beta \right\vert +\left\vert \gamma \right\vert \leq M}I_{\alpha \beta
\gamma }\left( \tau ,\tau _{0}\right)
\end{multline*}%
where 
\begin{multline*}
I_{\alpha \beta \gamma }\left( \tau ,\tau _{0}\right) \\
=\int \left\Vert c_{\gamma }\left( \cdot \right) b_{\beta }\left( \cdot
-y\right) \left( \mathcal{F}^{-1}a_{\alpha }\left( \left( 1-\tau \right)
\cdot +\tau y\right) -\mathcal{F}^{-1}a_{\alpha }\left( \left( 1-\tau
_{0}\right) \cdot +\tau _{0}y\right) \right) \right\Vert _{L^{2}\left(
X\right) }^{2}\limfunc{d}y
\end{multline*}%
with $a_{\alpha }\in \mathcal{S}\left( X^{\ast }\right) ,$ $b_{\beta }\in 
\mathcal{S}\left( X\right) $ and $c_{\gamma }\in \mathcal{S}^{\infty }\left(
X\right) $. Next we apply Lemma \ref{c4} to obtain 
\begin{equation*}
\lim_{\tau \rightarrow \tau _{0}}I_{\alpha \beta \gamma }\left( \tau ,\tau
_{0}\right) =0.
\end{equation*}

The proof of $(\func{ii})$ is similar to the proof of $(\func{i})$.
\end{proof}

For $\tau \in \mathcal{U}$ we consider the operator $K_{a,b}\left( \tau
\right) :\mathcal{S}\left( X\right) \rightarrow \mathcal{S}^{\ast }\left(
X\right) $ associated to the kernel $\mathcal{K}_{a,b}\left( \cdot ,\cdot
;\tau \right) $. Also, for $j=0,1$ and $\tau \in \mathcal{U}_{j}$ we
consider the Hilbert-Schmidt operator $K_{a,b,c}^{j}\left( \tau \right) \in 
\mathcal{B}_{2}\left( L^{2}\left( X\right) \right) $ associated to the
kernel $\mathcal{K}_{a,b,c}^{j}\left( \cdot ,\cdot ;\tau \right) $. Then 
\begin{eqnarray*}
K_{a,b,c}^{0}\left( \tau \right) &=&K_{a,b}\left( \tau \right) c\left(
Q\right) \left( 1-\triangle \right) ^{m/2},\quad \tau \in \mathcal{U}_{0}, \\
K_{a,b,c}^{1}\left( \tau \right) &=&\left( 1-\triangle \right) ^{m/2}c\left(
Q\right) K_{a,b}\left( \tau \right) ,\quad \tau \in \mathcal{U}_{1}.
\end{eqnarray*}

Let us show the first equality. For $f,g\in \mathcal{S}\left( X\right) $ and 
$\tau \in \mathcal{U}_{0}$ we have 
\begin{eqnarray*}
\left\langle g,K_{a,b}\left( \tau \right) c\left( Q\right) \left(
1-\triangle \right) ^{m/2}f\right\rangle &=&\left\langle g\otimes \overline{%
c\left( Q\right) \left( 1-\triangle \right) ^{m/2}f},\mathcal{K}_{a,b}\left(
\cdot ,\cdot ;\tau \right) \right\rangle \\
&=&\left\langle g\otimes \overline{c}\left( Q\right) \left( 1-\triangle
\right) ^{m/2}\overline{f},\mathcal{K}_{a,b}\left( \cdot ,\cdot ;\tau
\right) \right\rangle \\
&=&\left\langle g\otimes \overline{f},\left( 1-\triangle \right)
^{m/2}\left( c\left( \cdot \right) \mathcal{K}_{a,b}\left( x,\cdot ;\tau
\right) \right) \left( y\right) \right\rangle \\
&=&\left\langle g,K_{a,b,c}^{0}\left( \tau \right) f\right\rangle .
\end{eqnarray*}%
The second equality can be done in the same manner. For $f,g\in \mathcal{S}%
\left( X\right) $ and $\tau \in \mathcal{U}_{1}$ we have 
\begin{eqnarray*}
\left\langle g,\left( 1-\triangle \right) ^{m/2}c\left( Q\right)
K_{a,b}\left( \tau \right) f\right\rangle &=&\left\langle \overline{c}\left(
Q\right) \left( 1-\triangle \right) ^{m/2}g,K_{a,b}\left( \tau \right)
f\right\rangle \\
&=&\left\langle \left( \overline{c}\left( Q\right) \left( 1-\triangle
\right) ^{m/2}g\right) \otimes \overline{f},\mathcal{K}_{a,b}\left( \cdot
,\cdot ;\tau \right) \right\rangle \\
&=&\left\langle g\otimes \overline{f},\left( 1-\triangle \right)
^{m/2}\left( c\left( \cdot \right) \mathcal{K}_{a,b}\left( \cdot ,y;\tau
\right) \right) \left( x\right) \right\rangle \\
&=&\left\langle g,K_{a,b,c}^{1}\left( \tau \right) f\right\rangle .
\end{eqnarray*}

We recall that $\psi _{r}$ is the unique solution within $\mathcal{S}^{\ast
}(X)$ for 
\begin{equation*}
\left( 1-\triangle _{X}\right) ^{\frac{r}{2}}\psi _{r}=\delta ,
\end{equation*}%
where $r\in 
%TCIMACRO{\U{211d} }%
%BeginExpansion
\mathbb{R}
%EndExpansion
$ and $\delta $ is the delta function. Then $\mathcal{F}\left( \psi
_{r}\right) \in \mathcal{S}^{-r}\left( X^{\ast }\right) $.

Next we choose $c=\left\langle \cdot \right\rangle ^{s/2}$ so that 
\begin{eqnarray*}
K_{a,b}\left( \tau \right) &=&K_{a,b,c}^{0}\left( \tau \right) \left(
1-\triangle \right) ^{-m/2}\left\langle Q\right\rangle ^{-s/2},\quad \tau
\in \mathcal{U}_{0}, \\
K_{a,b}\left( \tau \right) &=&\left\langle Q\right\rangle ^{-s/2}\left(
1-\triangle \right) ^{-m/2}K_{a,b,c}^{1}\left( \tau \right) ,\quad \tau \in 
\mathcal{U}_{1}.
\end{eqnarray*}

We notice that $\left( \left( 1-\triangle \right) ^{-m/2}\left\langle
Q\right\rangle ^{-s/2}\right) ^{\ast }=\left\langle Q\right\rangle
^{-s/2}\left( 1-\triangle \right) ^{-m/2}$ and the kernel $\mathcal{K}_{s,m}$
of the operator $\left\langle Q\right\rangle ^{-s/2}\left( 1-\triangle
\right) ^{-m/2}$ is $\mathcal{K}_{s,m}\left( x,y\right) =\left\langle
x\right\rangle ^{-s/2}\psi _{m}\left( x-y\right) $. Since $s>n$, $m>n/2$ and 
$\mathcal{F}\left( \psi _{m}\right) \in \mathcal{S}^{-m}\left( X^{\ast
}\right) $ it follows feom Corollary \ref{c1} that $\mathcal{K}_{s,m}\in $ $%
L^{2}\left( X\times X\right) $ so that $\left\langle Q\right\rangle
^{-s/2}\left( 1-\triangle \right) ^{-m/2}$ and $\left( 1-\triangle \right)
^{-m/2}\left\langle Q\right\rangle ^{-s/2}$ are Hilbert-Schmidt operators.
From this we get the following extension of Cordes' lemma as a corollary.

\begin{corollary}[Cordes]
\label{c5}Let $s,t>n$, $a\in \mathcal{S}^{-t}\left( X^{\ast }\right) $, $%
b\in \mathcal{S}^{-s}\left( X\right) $ and 
\begin{equation*}
g:X\times X^{\ast }\rightarrow 
%TCIMACRO{\U{2102} }%
%BeginExpansion
\mathbb{C}
%EndExpansion
,\quad g\left( x,p\right) =\left( \mathcal{F}^{-1}a\right) \left( x\right)
\left( \mathcal{F}b\right) \left( p\right) ,\quad \left( x,p\right) \in
X\times X^{\ast },
\end{equation*}%
i.e. $g=\mathcal{F}^{-1}a\otimes \mathcal{F}b$. If $\tau \in \mathcal{U}$,
then $g_{X}^{\tau }\left( Q,P\right) $ has an extension in $\mathcal{B}%
_{1}\left( \mathcal{H}\left( X\right) \right) $ denoted also by $g_{X}^{\tau
}\left( Q,P\right) $. The mapping%
\begin{equation*}
\mathcal{U}\ni \tau \rightarrow g_{X}^{\tau }\left( Q,P\right) \in \mathcal{B%
}_{1}\left( \mathcal{H}\left( X\right) \right)
\end{equation*}%
is continuous.
\end{corollary}

\begin{proof}
It suffices to note that%
\begin{equation*}
\mathcal{K}_{g_{X}^{\tau }\left( Q,P\right) }=\left( \left( \limfunc{id}%
\otimes \mathcal{F}^{-1}\right) g\right) \circ C_{\tau }=\left( \mathcal{F}%
^{-1}a\otimes b\right) \circ C_{\tau }=\mathcal{K}_{a,b}\left( \cdot ,\cdot
;\tau \right) .
\end{equation*}
\end{proof}

\section{Schatten-class properties of pseudo-differential operators}

We are now able to consider $L^{2}$-boundedness and Schatten-class
properties of certain pseudo-differential operators. In the notation of \cite%
{Arsu}, $\mathfrak{S}$ will be $T^{\ast }\left( X\right) $ with the standard
symplectic structure and $\left( X,\left\vert \cdot \right\vert _{X}\right) $
an euclidean space, $(\mathcal{H}\left( X\right) ,\mathcal{W})$ will be the
Schr\"{o}dinger representation associated to the symplectic space $\mathfrak{%
S}$. Let $X=X_{1}\oplus ...\oplus X_{k}$ be an orthogonal decomposition and $%
X^{\ast }=X_{1}^{\ast }\oplus ...\oplus X_{k}^{\ast }$ be the dual
orthogonal decomposition. Let $\mathcal{U}_{0}$ be the set of all invetible
endomorphisms of $X$, $\mathcal{U}_{1}=1_{X}+\mathcal{U}_{0}$ and $\mathcal{U%
}=\mathcal{U}_{0}\cup \mathcal{U}_{1}$.

\begin{theorem}
\label{tcp2}Let $a\in \mathcal{S}^{\ast }(\mathfrak{S})$ and $1\leq p<\infty 
$. Assume that there are $\mathbf{t}=\left( t_{1},...,t_{k}\right) $, $%
\mathbf{s}=\left( s_{1},...,s_{k}\right) $ such that $t_{1},s_{1}>\frac{\dim
X_{1}}{4},...,t_{k},s_{k}>\frac{\dim X_{k}}{4}$ and 
\begin{equation*}
c=\left( 1-\triangle _{X_{1}}\right) ^{2t_{1}}\otimes ...\otimes \left(
1-\triangle _{X_{k}}\right) ^{2t_{k}}\otimes \left( 1-\triangle
_{X_{1}^{\ast }}\right) ^{2s_{1}}\otimes ...\otimes \left( 1-\triangle
_{X_{k}^{\ast }}\right) ^{2s_{k}}a\in L^{p}\left( \mathfrak{S}\right) .
\end{equation*}%
If $\tau \in \mathcal{U}$, then $a_{X}^{\tau }\left( Q,P\right) $ has an
extension in $\mathcal{B}_{p}\left( \mathcal{H}\left( X\right) \right) $
denoted also by $a_{X}^{\tau }\left( Q,P\right) $. The mapping%
\begin{equation*}
\mathcal{U}\ni \tau \rightarrow a_{X}^{\tau }\left( Q,P\right) \in \mathcal{B%
}_{p}\left( \mathcal{H}\left( X\right) \right)
\end{equation*}%
is continuous and for any $K$ a compact subset of $\mathcal{U}$, there is $%
C_{K}>0$ such that 
\begin{equation*}
\left\Vert a_{X}^{\tau }\left( Q,P\right) \right\Vert _{p}\leq
C_{K}\left\Vert c\right\Vert _{L^{p}\left( \mathfrak{S}\right) },
\end{equation*}%
for any $\tau \in K$.
\end{theorem}

\begin{proof}
If we use Corollary \ref{c5} insted of Corollary 5.3 in \cite{Arsu}, the
proof of this theorem is essentially the same as the proof of Theorem 6.1 in 
\cite{Arsu}.
\end{proof}

If we replace the $L^{p}$-conditions by $L^{\infty }$-conditions, then we
obtain the theorem on $L^{2}$-boundedness of Cordes' type.

\begin{theorem}
\label{tcp3}Let $a\in \mathcal{S}^{\ast }(\mathfrak{S})$. Assume that there
are $\mathbf{t}=\left( t_{1},...,t_{k}\right) $, $\mathbf{s}=\left(
s_{1},...,s_{k}\right) $ such that $t_{1},s_{1}>\frac{\dim X_{1}}{4}%
,...,t_{k},s_{k}>\frac{\dim X_{k}}{4}$ and 
\begin{equation*}
c=\left( 1-\triangle _{X_{1}}\right) ^{2t_{1}}\otimes ...\otimes \left(
1-\triangle _{X_{k}}\right) ^{2t_{k}}\otimes \left( 1-\triangle
_{X_{1}^{\ast }}\right) ^{2s_{1}}\otimes ...\otimes \left( 1-\triangle
_{X_{k}^{\ast }}\right) ^{2s_{k}}a\in L^{\infty }\left( \mathfrak{S}\right) .
\end{equation*}%
If $\tau \in \mathcal{U}$, then $a_{X}^{\tau }\left( Q,P\right) $ has an
extension in $\mathcal{B}\left( \mathcal{H}\left( X\right) \right) $ denoted
also by $a_{X}^{\tau }\left( Q,P\right) $. The mapping%
\begin{equation*}
%TCIMACRO{\U{211d} }%
%BeginExpansion
\mathbb{R}
%EndExpansion
\ni \tau \rightarrow a_{X}^{\tau }\left( Q,P\right) \in \mathcal{B}\left( 
\mathcal{H}\left( X\right) \right)
\end{equation*}%
is continuous and for any $K$ a compact subset of $\mathcal{U}$, there is $%
C_{K}>0$ such that 
\begin{equation*}
\left\Vert a_{X}^{\tau }\left( Q,P\right) \right\Vert _{\mathcal{B}\left( 
\mathcal{H}\left( X\right) \right) }\leq C_{K}\left\Vert c\right\Vert
_{L^{\infty }\left( \mathfrak{S}\right) },
\end{equation*}%
for any $\tau \in K$.
\end{theorem}

The proof of this theorem is essentially the same, the only difference being
the reference to the new $\tau $-version of Cordes' lemma $\ $insted of the
original one.

We recall some notations from \cite{Arsu}. We consider the symplectic space $%
\mathfrak{S}=T^{\ast }\left( X\right) $. An orthogonal decomposition of $X$, 
$X=X_{1}\oplus ...\oplus X_{k}$, gives an orthogonal decomposition of $%
\mathfrak{S}$, $\mathfrak{S}=X_{1}\oplus ...\oplus X_{k}\oplus X_{1}^{\ast
}\oplus ...\oplus X_{k}^{\ast }$, if on $\mathfrak{S}$ we consider the
euclidean norm $\left\Vert \left( x,p\right) \right\Vert _{\mathfrak{S}%
}^{2}=\left\Vert x\right\Vert _{X}^{2}+\left\Vert p\right\Vert _{X^{\ast
}}^{2}$. We shall choose an orthonormal basis in each space $X_{j}$, $%
j=1,...,k$, while in $X_{j}^{\ast }$, $j=1,...,k$ we shall consider the dual
bases. Then $\partial ^{X}=\left( \partial ^{X_{1}},...,\partial
^{X_{k}}\right) $, $\partial ^{X^{\ast }}=\left( \partial ^{X_{1}^{\ast
}},...,\partial ^{X_{k}^{\ast }}\right) .$

For $1\leq p\leq \infty $ and $\mathbf{t}=\left( t_{1},...,t_{k}\right) \in 
%TCIMACRO{\U{2115} }%
%BeginExpansion
\mathbb{N}
%EndExpansion
^{k}$ we set $\mathcal{M}_{\mathbf{t}}^{p}=\mathcal{M}_{t_{1},...,t_{k}}^{p}$
for the space of all distributions $a\in \mathcal{S}^{\ast }(\mathfrak{S})$
whose derivatives $\left( \partial ^{X_{1}}\right) ^{\alpha _{1}}...\left(
\partial ^{X_{k}}\right) ^{\alpha _{k}}\left( \partial ^{X_{1}^{\ast
}}\right) ^{\beta _{1}}...\left( \partial ^{X_{k}^{\ast }}\right) ^{\beta
_{k}}a$ belong to $L^{p}\left( \mathfrak{S}\right) $ when $\alpha _{j},\beta
_{j}\in 
%TCIMACRO{\U{2115} }%
%BeginExpansion
\mathbb{N}
%EndExpansion
^{\dim X_{j}}$, $\left\vert \alpha _{j}\right\vert ,\left\vert \beta
_{j}\right\vert \leq t_{j}$, $j=1,...,k$. On this space we shall consider
the natural norm $\left\vert \cdot \right\vert _{p,\mathbf{t}}=\left\vert
\cdot \right\vert _{p,t_{1},...,t_{k}}$ defined by 
\begin{equation*}
\left\vert a\right\vert _{p,\mathbf{t}}=\max_{\left\vert \alpha
_{1}\right\vert ,\left\vert \beta _{1}\right\vert \leq t_{1},...,\left\vert
\alpha _{k}\right\vert ,\left\vert \beta _{k}\right\vert \leq
t_{k}}\left\Vert \left( \partial ^{X_{1}}\right) ^{\alpha _{1}}...\left(
\partial ^{X_{k}}\right) ^{\alpha _{k}}\left( \partial ^{X_{1}^{\ast
}}\right) ^{\beta _{1}}...\left( \partial ^{X_{k}^{\ast }}\right) ^{\beta
_{k}}a\right\Vert _{L^{p}}.
\end{equation*}

Let $m_{1}=\left[ \frac{\dim X_{1}}{2}\right] +1,...,m_{k}=\left[ \frac{\dim
X_{k}}{2}\right] +1$ and $\mathbf{m=}\left( m_{1},...,m_{k}\right) $.

A consequence of Theorem \ref{tcp2} and of Lemma 6.3 in \cite{Arsu} is the
following

\begin{theorem}
Assume that $1\leq p<\infty $ and let $a\in \mathcal{S}^{\ast }(\mathfrak{S}%
) $. If $a\in \mathcal{M}_{2m_{1},...,2m_{k}}^{p}$, then for any $\tau \in 
\mathcal{U}$, $a_{X}^{\tau }\left( Q,P\right) $ has an extension in $%
\mathcal{B}_{p}\left( \mathcal{H}\left( X\right) \right) $ denoted also by $%
a_{X}^{\tau }\left( Q,P\right) $. The mapping 
\begin{equation*}
\mathcal{U}\ni \tau \rightarrow a_{X}^{\tau }\left( Q,P\right) \in \mathcal{B%
}_{p}\left( \mathcal{H}\left( X\right) \right)
\end{equation*}%
is continuous and for any $K$ a compact subset of $\mathcal{U}$, there is $%
C_{K}>0$ such that 
\begin{equation*}
\left\Vert a_{X}^{\tau }\left( Q,P\right) \right\Vert _{p}\leq
C_{K}\left\vert a\right\vert _{p,2m_{1},...,2m_{k}},
\end{equation*}%
for any $\tau \in K$.
\end{theorem}

Similarly, for $p=\infty $, a consequence of Theorem \ref{tcp3} and of Lemma
6.3 in \cite{Arsu} is the celebrated Calderon-Vaillancourt Theorem.

\begin{theorem}[Calderon, Vaillancourt]
If $a\in \mathcal{M}_{2m_{1},...,2m_{k}}^{\infty }$, then $a_{X}^{\tau
}\left( Q,P\right) $ is $L^{2}$-bounded for any $\tau \in \mathcal{U}$. The
mapping 
\begin{equation*}
\mathcal{U}\ni \tau \rightarrow a_{X}^{\tau }\left( Q,P\right) \in \mathcal{B%
}\left( \mathcal{H}\left( X\right) \right)
\end{equation*}%
is continuous and for any $K$ a compact subset of $\mathcal{U}$, there is $%
C_{K}>0$ such that 
\begin{equation*}
\left\Vert a_{X}^{\tau }\left( Q,P\right) \right\Vert _{\mathcal{B}\left( 
\mathcal{H}\left( X\right) \right) }\leq C_{K}\left\vert a\right\vert
_{\infty ,2m_{1},...,2m_{k}},
\end{equation*}%
for any $\tau \in K$.
\end{theorem}

The next two theorems are consequences of Lemma 6.6 in \cite{Arsu}, Theorem %
\ref{tcp2} and Theorem \ref{tcp3}. Recall that the Sobolev space $%
H_{p}^{s}\left( \mathfrak{S}\right) $, $s\in 
%TCIMACRO{\U{211d} }%
%BeginExpansion
\mathbb{R}
%EndExpansion
$, $1\leq p\leq \infty $, consists of all $a\in \mathcal{S}^{\ast }(%
\mathfrak{S})$ such that $\left( 1-\triangle _{\mathfrak{S}}\right)
^{s/2}a\in L^{p}\left( \mathfrak{S}\right) $, and we set $\left\Vert
a\right\Vert _{H_{p}^{s}\left( \mathfrak{S}\right) }\equiv \left\Vert \left(
1-\triangle _{\mathfrak{S}}\right) ^{s/2}a\right\Vert _{L^{p}\left( 
\mathfrak{S}\right) }$.

\begin{theorem}
Assume that $1\leq p<\infty $. If $s>2\dim X$ and $a\in H_{p}^{s}\left( 
\mathfrak{S}\right) $, then for any $\tau \in \mathcal{U}$, $a_{X}^{\tau
}\left( Q,P\right) $ has an extension in $\mathcal{B}_{p}\left( \mathcal{H}%
\left( X\right) \right) $ denoted also by $a_{X}^{\tau }\left( Q,P\right) $.
The mapping 
\begin{equation*}
\mathcal{U}\ni \tau \rightarrow a_{X}^{\tau }\left( Q,P\right) \in \mathcal{B%
}_{p}\left( \mathcal{H}\left( X\right) \right)
\end{equation*}%
is continuous and for any $K$ a compact subset of $\mathcal{U}$, there is $%
C_{K}>0$ such that 
\begin{equation*}
\left\Vert a_{X}^{\tau }\left( Q,P\right) \right\Vert _{p}\leq
C_{K}\left\Vert a\right\Vert _{H_{p}^{s}\left( \mathfrak{S}\right) },
\end{equation*}%
for any $\tau \in K$.
\end{theorem}

\begin{theorem}
If $s>2\dim X$ and $a\in H_{\infty }^{s}\left( \mathfrak{S}\right) $, then $%
a_{X}^{\tau }\left( Q,P\right) $ is $L^{2}$-bounded for any $\tau \in 
\mathcal{U}$. The mapping 
\begin{equation*}
\mathcal{U}\ni \tau \rightarrow a_{X}^{\tau }\left( Q,P\right) \in \mathcal{B%
}\left( \mathcal{H}\left( X\right) \right)
\end{equation*}%
is continuous and for any $K$ a compact subset of $\mathcal{U}$, there is $%
C_{K}>0$ such that 
\begin{equation*}
\left\Vert a_{X}^{\tau }\left( Q,P\right) \right\Vert _{\mathcal{B}\left( 
\mathcal{H}\left( X\right) \right) }\leq C_{K}\left\Vert a\right\Vert
_{H_{\infty }^{s}\left( \mathfrak{S}\right) },
\end{equation*}%
for any $\tau \in K$.
\end{theorem}

If we note that $a_{X}^{\tau }\left( Q,P\right) \in \mathcal{B}_{2}\left( 
\mathcal{H}\left( X\right) \right) $ whenever $a\in L^{2}\left( \mathfrak{S}%
\right) =H_{2}^{0}\left( \mathfrak{S}\right) $, then the last two theorems
and standard interpolation results in Sobolev spaces (see \cite[Theorem 6.4.5%
]{Bergh}) give us the following

\begin{theorem}
Let $\mu >1$, $1\leq p<\infty $ and $n=\dim X$ . If $a\in H_{p}^{2\mu
n\left\vert 1-2/p\right\vert }\left( \mathfrak{S}\right) $, then for any $%
\tau \in \mathcal{U}$, $a_{X}^{\tau }\left( Q,P\right) $ has an extension in 
$\mathcal{B}_{p}\left( \mathcal{H}\left( X\right) \right) $ denoted also by $%
a_{X}^{\tau }\left( Q,P\right) $. The mapping 
\begin{equation*}
\mathcal{U}\ni \tau \rightarrow a_{X}^{\tau }\left( Q,P\right) \in \mathcal{B%
}_{p}\left( \mathcal{H}\left( X\right) \right)
\end{equation*}%
is continuous and for any $K$ a compact subset of $\mathcal{U}$, there is $%
C_{K}>0$ such that 
\begin{equation*}
\left\Vert a_{X}^{\tau }\left( Q,P\right) \right\Vert _{p}\leq
C_{K}\left\Vert a\right\Vert _{H_{p}^{2\mu n\left\vert 1-2/p\right\vert
}\left( \mathfrak{S}\right) },
\end{equation*}%
for any $\tau \in K$.
\end{theorem}


\begin{thebibliography}{99}
\bibitem{Arsu} G. Arsu: On Schatten-von Neumann class properties of
pseudo-differential operators. The Cordes-Kato method, preprint at
http://arxiv.org/abs/math.AP/0511080. To appear in Journal of Operator
Theory.

\bibitem{Beals} R. Beals: On the boundedness of pseudodifferential
operators, Comm. Partial Differential Equations 2 (1977), 1063- 1070.

\bibitem{Bergh} J. Bergh and J. L\"{o}fstr\"{o}m: Interpolation Spaces,
Springer, Berlin, 1976.

\bibitem{Birman} M.S. Birman and M.Z. Solomjak: Spectral Theory of
Self-Adjoint Operators in HilbertSpace, D. Riedel Publishing Company, 1987.

\bibitem{Boulkhemair} A. Boulkhemair: $L^{2}$ Estimates for Weyl
Quantization, J. Funct. Anal. 165 (1999), 173-204.

\bibitem{Bourdaud} G. Bourdaud and Y. Meyer: In\'{e}galit\'{e}s $L^{2}$ pr%
\'{e}cis\'{e}es pour la classe $S_{0,0}^{0}$, Bull. Soc. Math. France 116
(1988), 401-412.

\bibitem{Georgescu} A. Boutet de Monvel-Berthier and V. Georgescu: Graded
C*-algebras associated to symplectic spaces and spectral analysis of many
channel Hamiltonians, Dynamics of complex and irregular systems (Bielefeld,
1991), 22-66, Bielefeld Encount. Math. Phys., VIII, World Sci. Publishing,
1993.

\bibitem{Calderon} A.P. Calder\'{o}n and R. Vaillancourt: On the boundedness
of pseudo-differential operators, J. Math. Soc. Japan 23 (1971)\ 374-378.

\bibitem{Childs} A.G. Childs: On $L^{2}$ boundedness of pseudo-differential
operators, Proc. Amer. Math. Soc. 61 (1976), 252-254.

\bibitem{Coifman} R. Coifman and Y. Meyer: Au del\`{a} des op\'{e}rateurs
pseudo-diff\'{e}rentiels, Ast\'{e}risque 57 (1978).

\bibitem{Cordes} H.O. Cordes: On compactness of commutators of
multiplications and convolutions, and boundedness of pseudodifferential
operators, J. Funct. Anal. 18 (1975), 115-131.

\bibitem{Folland} G. B. Folland: Harmonic analysis in phase space, Princeton
University Press, 1989.

\bibitem{Grochenig} K. Gr\"{o}chenig and C. Heil: Modulation spaces and
pseudodifferential operators, Integral Equations Operator Theory 34 (1999),
439-457.

\bibitem{Heil} C. Heil, J. Ramanathan, and P. Topiwala: Singular values of
compact pseudodifferential operators, J. Funct. Anal. 150 (1997), 426-452.

\bibitem{Hormander1} L. H\"{o}rmander: The Weyl calculusof
pseudo-differential operators, Comm. Pure Appl. Math. 32 (1979), 359-443.

\bibitem{Hormander2} L. H\"{o}rmander: The Analysis of Linear Partial
Differential Operators, vol. I, III, Springer-Verlag, Berlin Heidelberg New
York Tokyo, 1983, 1985.

\bibitem{Howe} R. Howe: Quantum mechanics and partial differential
equations, J. Funct. Anal. 38 (1980), 188-254.

\bibitem{Hwang} I.L. Hwang: On $L^{2}$-boundedness of pseudodifferential
operators, Trans. Amer. Math. Soc. 302, No.1 (1987), 55-76.

\bibitem{Kastler} D. Kastler: The C*-algebras of a free boson field. I.
Discussion of basic facts, Comm. Math. Phys. 9 (1965), 14-48.

\bibitem{Kato} T. Kato: Boundedness of some pseudo-differential operators,
Osaka J. Math. 13 (1976), 1-9.

\bibitem{Raeburn} I. Raeburn and D.P. Wiliams: Morita Equivalence and
Continuous Trace $C^{\ast }$-Algebras, American Mathematical Society, 1998.

\bibitem{Simon} M. Reed and Methods of Modern Mathematical Physics, vol. I,
Academic Press, New York San Francisco London, 1972.

\bibitem{Simon1} B. Simon: Trace Ideals and their Applications, vol.35,
London Mathematical Society Lecture Note Series, Cambrige University Press,
Cambrige London New York Melbourne, 1979.

\bibitem{Taylor} M.E. Taylor: Pseudodifferential Operators, Princeton,
University Press, 1981.

\bibitem{Triebel} H. Triebel: Interpolation Theory, Function Spaces,
Differential Operators, VEB Deutscher Verlag der Wissenschaften, Berlin,
1978.

\bibitem{Toft1} J. Toft: Continuity and Positivity Problems in
Pseudo-Differential Calculus, Thesis, Departament of Mathematics, University
of Lund, Lund, 1996.

\bibitem{Toft2} J. Toft: Regularizations, decompositions and lower bound
problems in the Weyl calculus, Comm. Partial Differential Equations 7\&8
(2000), 1201- 1234.

\bibitem{Toft3} J. Toft: Subalgebras to a Wienner type algebra of
pseudo-differential operators, Ann. Inst. Fourier 51 (2001), 1347-1383.

\bibitem{Toft4} J. Toft: Continuity properties in non-commutative
convolution algebras, with applications in pseudo-differential calculus,
Bull. Sci. Math. 126 (2002), 115-142.
\end{thebibliography}
\end{document}